\documentclass[12pt,a4paper]{article}   

\input amssym
 
\catcode`\@=11

\font\decp=cmr8

\font\tenmsa=msam10 \font\sevenmsa=msam7 \font\fivemsa=msam5
\font\tenmsb=msbm10 at 12pt\font\sevenmsb=msbm7 at 9pt\font\fivemsb=msbm5 at 7pt
\newfam\msafam \newfam\msbfam
\textfont\msafam=\tenmsa  \scriptfont\msafam=\sevenmsa \scriptscriptfont\msafam=
\fivemsa
\textfont\msbfam=\tenmsb  \scriptfont\msbfam=\sevenmsb \scriptscriptfont\msbfam=
\fivemsb

\font\teneufm=eufm10 \font\seveneufm=eufm7 \font\fiveeufm=eufm5
\newfam\eufmfam 
\textfont\eufmfam=\teneufm \scriptfont\eufmfam=\seveneufm\scriptscriptfont
\eufmfam=\fiveeufm

\def\hexnumber@#1{\ifcase#1 0\or1\or2\or3\or4\or5\or6\or7\or8\or9\or
 A\or B\or C\or D\or E\or F\fi }

\edef\msa@{\hexnumber@\msafam} \edef\msb@{\hexnumber@\msbfam}

\mathchardef\qed="0\msa@03

\def\Bbb#1{{\fam\msbfam\relax#1}}

\def\R{{\Bbb R}}    

\catcode`\@=12

\newtheorem{teo}{Theorem}[section]
\newtheorem{pro}[teo]{Proposition}
\newtheorem{cor}[teo]{Corollary} 
\newtheorem{eje}[teo]{Examples} 
\newtheorem{lem}[teo]{Lemma} 
\newtheorem{obs}[teo]{Remark} 
\newtheorem{defin}[teo]{Definition}

\newcommand{\bteo}{\begin{teo}}
\newcommand{\eteo}{\end{teo}}
\newcommand{\bpro}{\begin{pro}}
\newcommand{\epro}{\end{pro}}
\newcommand{\bcor}{\begin{cor}}
\newcommand{\ecor}{\end{cor}}
\newcommand{\blem}{\begin{lem}}
\newcommand{\elem}{\end{lem}}
\newcommand{\bdefi}{\begin{defin}\rm}
\newcommand{\edefi}{\end{defin}}
\newcommand{\beje}{\begin{eje}\rm}
\newcommand{\eeje}{\end{eje}}
\newcommand{\bobs}{\begin{obs}\rm}
\newcommand{\eobs}{\end{obs}}
\newcommand{\bequ}{\begin{equation}}
\newcommand{\eequ}{\end{equation}}

\def\bdem{{\sl Proof. }}
\def\edem{\par\bigskip}
\def\lpw{{L^p(w)}}
\def\lopw{{\Lambda^p(w)}}
\def\loqw{{\Lambda^q(w)}}
\def\louw{{\Lambda^1(w)}}
\def\gamuiw{{\Gamma^{1,\infty}(w)}}

\def\gamuw{{\Gamma^{1}(w)}}
\def\gamuv{{\Gamma^{1}(v)}}
\def\gamupw{{\Gamma^{1,p}(w)}}
\def\gamuqw{{\Gamma^{1,q}(w)}}
\def\gampqw{{\Gamma^{p,q}(w)}}
\def\gampw{{\Gamma^{p}(w)}}

\def\gamqv{{\Gamma^{q}(v)}}
\def\gamqw{{\Gamma^{q}(w)}}
\def\gampiw{{\Gamma^{p,\infty}(w)}}
\def\gamqiw{{\Gamma^{q,\infty}(w)}}
\def\lopiw{{\Lambda^{p,\infty}(w)}}
\def\loqiw{{\Lambda^{q,\infty}(w)}}
\def\lpwd{L^p_{\hbox{\decp dec}}(w)}
\def\luwd{L^1_{\hbox{\decp dec}}(w)}

\def\luiw{L^{1,\infty}(w)}
\def\lpiw{L^{p,\infty}(w)}
\def\gamqaw{{\Gamma^{q}_{\alpha}(w)}}
\def\gampaw{{\Gamma^{p}_{\alpha}(w)}}

\def\gamquw{{\Gamma^{q}_{1}(w)}}
\def\gamuaw{{\Gamma^{1}_{\alpha}(w)}}

\begin{document}
\title{New Lorentz spaces for the restricted weak-type Hardy's inequalities}  
\author{ {\bf Joaquim Mart\'\i n$^*$} \\Dept.\
Appl.\ Math.\  and Analysis\\ University of  Barcelona\\
E-08071 Barcelona, SPAIN\\E-mail: {\tt jmartin@mat.ub.es}   \and {\bf Javier Soria\thanks{Partially supported
 by DGICYT PB97-0986 and CIRIT 1999SGR00061} }\\Dept.\
Appl.\ Math.\  and Analysis\\ University of  Barcelona\\
E-08071 Barcelona, SPAIN\\E-mail: {\tt soria@mat.ub.es}}
\date{} 
\maketitle
\begin{center}{\bf Abstract}
\end{center}
Associated to the class of restricted-weak type weights for the Hardy operator $R_p$, we find a new class of Lorentz spaces
for which the normability property holds. This result is analogous to the  characterization given by Sawyer for the
classical  Lorentz spaces. We also show that these new spaces are very natural to study the existence of equivalent norms
described in terms of the maximal function.
\vskip 1cm
\noindent
{\bf Mathematics Subject Classification 2000:} 46E30, 26D10

\vskip 1cm
\noindent
{\bf Keywords:} Hardy operator, Lorentz spaces, monotone functions, weighted inequalities.
\newpage

\section{Introduction}

During the last decade, the characterization of the normability of all sorts of weighted Lorentz spaces in terms of either
weak-type or strong inequalities for the Hardy operator, has been completely settled down (see \cite{AM}, \cite{Sa},
\cite{CGS}, and \cite{So}). One of the main results of this paper is to consider the class of weights $R_p$ introduced by
Neugebauer (\cite{Ne}), which characterizes the restricted weak-type boundedness of the Hardy operator, and show that there
exists a new satisfactory Lorentz space, in the sense that the normability condition is described in terms of $R_p$ (thus
completing the picture of equivalences between Banach Lorentz-type spaces and weighted inequalities). We will also prove that
these new spaces are useful to consider the problem of finding equivalent norms in the classical Lorentz spaces depending upon
the maximal function (the case $p=1$ was an open problem, which we now solve in its full generality).
\medskip

In section 2 we introduce the spaces $\gamqaw$, study the embeddings with respect to both $\loqw$ and $\gamqw$,  and
characterize the existence of an equivalent norm. In section 3 we prove that all Banach  Lorentz spaces $\louw$ admit a norm
depending on the maximal function, and give an explicit formula for it (see Theorem~\ref{gamuv}).  This study leads us to
investigate how the space
$\louw$ fits among  the range of Banach spaces
$\gamupw$, which arise in a natural way in the theory. As a consequence we find that,  on this scale,  there are not
intermediate  weighted inequalities for the Hardy operator, that is, as soon as one gets  a better estimate than the weak-type
inequality, we obtain the best possible estimate, namely the strong-type (1,1) inequality,  (see Theorem~\ref{gamup}).  We
finally give some examples and applications to show that we can have all sort of different situations for the embeddings
$\gamuw\subset\louw\subset
\gamuiw$.
\medskip

In what follows, we will use the notation decreasing (increasing) with the meaning nonincreasing (resp.\ nondecreasing).  We
will say that a function is positive whenever it is nonnegative. A weight is a positive measurable function, locally
integrable on $\R^+$. Given a weight $w$, we denote by $W(t)=\int_0^tw(s)\,ds$.  $f^*$ denotes the decreasing
rearrangement of the function
$f$.  $L^p_{\hbox{\decp dec}}(w)$ is the cone of positive and decreasing functions in $L^p(w)$ (similarly for
$L^{p,\infty}_{\hbox{\decp dec}}(w)$). Constants such as
$C$ may have different values from one occurrence to the next, but they will always be irrelevant for the arguments used. Two 
positive quantities $A$ and $B$, are said to be equivalent
($A\approx B$) if there exists a constant $C>1$ (independent of the essential parameters defining $A$ and
$B$) such that $C^{-1}A\le B\le CA$. For
other standard notations we refer the reader to \cite{BS}. 

\section{The spaces $\gamqaw$}

We begin by introducing the definitions of the spaces we are going to study, and give some well-known results which will be
useful in our proofs.

\bdefi\label{lore}
Let $0<p<\infty$  and $w$ be a weight. Then we define the weighted Lorentz space
$$
\lopw=\{f:\R^n\longrightarrow \R^+; \Vert f\Vert_{\lopw}<\infty\},
$$
where,
\bequ\label{norm}
\Vert f\Vert_{\lopw}=\bigg(\int_0^{\infty}(f^*(t))^pw(t)\,dt\bigg)^{1/p}.
\eequ
\edefi

Classical examples are obtained when one considers power weights. Thus, if $w(t)=t^{(p/q)-1}$ then $\lopw=L^{q,p}$, which is
the Lebesgue space $L^p$ if $p=q$.  In general  (\ref{norm}) does not define a norm. In fact, Lorentz proved (\cite{Lo}) that 
(\ref{norm}) is a norm, if and only if $p\ge 1$ and $w$ is a decreasing function. This result was later on improved by Sawyer
(\cite{Sa}) by giving a characterization of the normability of $\lopw$, $p>1$. In order to formulate this result, we need first to
introduce some important operators.

\bdefi
The Hardy operator is defined  as
\bequ\label{hardy}
Sf(t)=\displaystyle{1\over t}\int_0^tf(s)\,ds,\qquad t>0.
\eequ
The maximal function of $f$ is
\bequ\label{maxfun}
f^{**}(t)=Sf^*(t).
\eequ
\edefi

The main result about the boundedness of the Hardy operator that we need, is the following theorem due to Ari\~no and 
Muckenhoupt:

\bteo\label{bpam}
For $0<p<\infty$, the following conditions are equivalent:
\medskip

(i) $S:\lpwd\longrightarrow L^p(w)$,

\medskip
(ii) $w\in B_p$; i.e., for all $r>0$,
\bequ\label{bp}
\int_r^{\infty}{w(s)\over s^p}\,ds\le {C\over r^p}\int_0^r w(s)\,ds.
\eequ
\eteo 

There are also some
other Lorentz-type spaces that we are going to consider:

\bdefi\label{gamupw}
Let $0<p<\infty$,  $0<q\le\infty$, and $w$ be a weight. Then we define:

\medskip\noindent
(i) The weak-type weighted Lorentz space,
$$
\lopiw=\{f:\R^n\longrightarrow \R^+; \Vert f\Vert_{\lopiw}<\infty\},
$$
where,
\bequ\label{normd}
\Vert f\Vert_{\lopiw}=\sup_{t>0}(f^*(t)W^{1/p}(t)).
\eequ

\medskip\noindent
(ii) The Gamma space,
$$
\gampqw=\{f:\R^n\longrightarrow \R^+; \Vert f\Vert_{\gampqw}<\infty\},
$$
where,
\bequ\label{normg}
\Vert f\Vert_{\gampqw}=\bigg(\int_0^{\infty}(f^{**}(t))^q(W(t))^{(q/p)-1}w(t)\,dt\bigg)^{1/q},
\eequ
if $q<\infty$, and if $q=\infty$, 
\bequ\label{normgd}
\Vert f\Vert_{\gampiw}=\sup_{t>0}(f^{**}(t)W^{1/p}(t)).
\eequ
In case $p=q$ we simply write $\gampw$.
\edefi
\noindent
The importance of (\ref{maxfun}) is the subadditivity property $(f+g)^{**}(t)\le f^{**}(t)+g^{**}(t)$.  This immediately
implies that (\ref{normg})  defines a (complete) norm if $p\ge 1$ (in  (\ref{normgd}) this is true for $p>0$), and one always
has that
$\gampw\subset\lopw\subset\lopiw$ (see  \cite{CPSS} for more information on this kind of embeddings).  We can now write
Sawyer's theorem:

\bteo\label{normlop}
If $1<p$, then the following conditions are equivalent:
\medskip

(i) $\lopw$ is a Banach space.

\medskip
(ii) $w\in B_p$.

\medskip
(iii) $\lopw=\gampw$.
\eteo

Thus, as a consequence we have that if $\lopw$ is normable, then the equivalent norm is always described in terms of the
maximal function (in fact $\Vert f\Vert_{\lopw}\approx\Vert f\Vert_{\gampw}).$  Similarly, the normability of the
weak-type spaces $\lopiw$ was also characterized in \cite{So}:

\bteo\label{normlopi}
If $0<p$, then the following conditions are equivalent:
\medskip

(i) $\lopiw$ is a Banach space.

\medskip
(ii) $w\in B_p$.

\medskip
(iii) $\lopiw=\gampiw$.
\eteo

Thus, as in Theorem~\ref{normlop}, we also have in this case that the equivalent norm in $\lopiw$ is given in terms of the
maximal function. In view of these results, we investigate what happens in the only case that is left, namely $\louw$. For this
endpoint space, we recall  the normability characterization proved in \cite{CGS}:

\bteo\label{normlouw}
The following conditions are equivalent:
\medskip
 
(i) $\louw$ is a Banach space.

\medskip
(ii)  $\louw\subset\gamuiw$.

\medskip
(iii) $\displaystyle{{1\over t}\int_0^tw(r)\,dr\le{C\over s}\int_0^s w(r)\,dr}$,  if  $0<s\le t$.
\eteo

\bobs\label{debilfuerte}
As a corollary of this theorem we have that if $\louw$ is normable, there exists a decreasing weight $\tilde w$ such that
$\louw=\Lambda^1(\tilde w)$, and hence, in many cases we can assume without loss of generality, that $w$ is already a
decreasing weight.  Also, it can be shown that (iii) of the previous theorem is equivalent to the weak-type inequality
$S:\luwd\longrightarrow\luiw$, and also that the $B_p$ condition is equivalent to $S:\lpwd\longrightarrow\lpiw$, if $p>1$.
Hence the above theorems can be summarized as follows:

\medskip
\noindent
- $\lopw$ is a Banach space, if and only if $p\ge 1$ and $S:\lpwd\to\lpiw$.

\medskip
\noindent
- $\lopiw$ is a Banach space, if and only if  $S:\lpwd\to\lpw$.

\medskip
There exists a third class of inequalities, considered in \cite{Ne}, which are referred as $R_p$, where the restricted 
weak-type boundedness for the Hardy operator and monotone functions (i.e., $w(\{S\chi_{(0,r)}>\lambda\})\le
CW(r)/\lambda^p$), is characterized by the condition

\bequ\label{rp}
{1\over s^p}\int_0^sw(x)\,dx\le {C\over r^p}\int_0^rw(x)\,dx,
\eequ
for $0<r<s<\infty$.

\eobs

We will show that there exist a suitable class of Lorentz spaces for which (\ref{rp}) is the necessary and sufficient condition for
the normability.

\bdefi\label{loalfap}
Let $0<p<\infty$, $0\le\alpha\le p$,  and $w$ be a weight. Then we define the  space
$$
\gampaw=\{f:\R^n\longrightarrow \R^+; \Vert f\Vert_{\gampaw}<\infty\},
$$
where,
\bequ\label{normalfa}
\Vert f\Vert_{\gampaw}=\bigg(\int_0^{\infty}(f^*(t))^{\alpha}(f^{**}(t))^{p-\alpha}w(t)\,dt\bigg)^{1/p}.
\eequ
\edefi

\bobs\label{motiv}
Observe that $\Gamma^p_0(w)=\gampw$,  $\Gamma^p_p(w)=\lopw$, and
$\gampw\subset\gampaw\subset\Gamma^p_{\beta}(w)\subset\lopw$, if $0\le\alpha\le\beta\le p$. It is also easy to prove
that if $1\le p<\infty$, and $0<\alpha\le p$, then $\gampw=\gampaw$, if and only if $w\in B_p$.
\eobs

We are going to consider first  $\alpha=1$. For this case we need some previous results which are of independent
interest.

\blem\label{wq} If $w$ is a decreasing weight such that $w(\infty)=0$, then for every $1\le q<\infty$, there exists a weight
$w_q$ such that,
\bequ\label{ecuwq}
\bigg(\int_0^rw(x)\,dx\bigg)^q\approx\int_0^rw_q(x)\,dx+r^q\int_r^{\infty}w_q(x)\,{dx\over x^q}.
\eequ
\elem

\bdem
We are going to assume first that $w\in {\cal C}^1(\R^+)$. Then, taking 
\bequ\label{defwq}
w_q(r)=-r^q{{\rm d}\over
{\rm dx}}\bigg({W^{q-1}w\over x^{q-1}}\bigg)(r),
\eequ
which defines a positive function, we have that
$$
\int_r^{\infty}{w_q(s)\over s^q}\,ds={W^{q-1}(r)w(r)\over r^{q-1}},
$$
and, since $0\le\lim_{s\to 0}sW^{q-1}(s)w(s)\le\lim_{s\to 0}W^q(s)=0$, then
$$
\int_0^r w_q(s)\,ds=-rW^{q-1}(r)w(r)+\bigg(\int_0^rw(s)\,ds\bigg)^q,
$$
which proves (\ref{ecuwq}). If $w$ is only a continuous function we define the auxiliary function
$$
\Phi(t)={1\over t}\int_t^{2t}W(s)\,ds=\int_1^{2}W(ts)\,ds.
$$
It is clear that $\Phi\in{\cal C}^2(\R^+)$, and we only need to show that $\Phi$ is concave, $\Phi\approx W$, and
$\lim_{t\to 0}t\Phi'(t)=\lim_{t\to\infty }\Phi'(t)=0$:

Take $a,b\in[0,1]$, $a+b=1$, and $s,t>0$. Then 
$$
a\Phi(s)+b\Phi(t)=\int_1^2(aW(sr)+bW(tr))\,dr\le\int_1^2W(asr+btr)\,dr=\Phi(as+bt),
$$
which shows the concavity. Also,
$$
W(t)\le\Phi(t)\le W(2t)\le 2W(t).
$$
Now,  since
$$
\Phi'(t)={2W(2t)-W(t)\over t}-{1\over t^2}\int_t^{2t}W(s)\,ds,
$$
then $\lim_{t\to 0}t\Phi'(t)=0$, and $\lim_{t\to\infty }\Phi'(t)=0$, since, by hypothesis,
$$
\lim_{t\to\infty}{1\over t^2}\int_t^{2t}W(s)\,ds=\lim_{t\to\infty}{\Phi(t)\over t}\approx \lim_{t\to\infty}{W(t)\over
t}=0.
$$
If we now define the function $w_q$ as in (\ref{defwq}), changing $W$ by $\Phi$,  we have that $w_q$ satisfies
(\ref{ecuwq}).  For a general weight
$w$, we repeat the previous argument twice. $\qquad\qed$
\edem

\bpro\label{teo1} Let $w$ be a decreasing weight.

\medskip
\noindent
(a) If $1<q<\infty$ and $v$ is a weight for which $\louw\subset\gamqv\subset\gamuiw$, then 
$$
\int_0^{\infty}(f^{**}(x))^qv(x)\,dx\approx\int_0^{\infty}f^*(x)(f^{**}(x))^{q-1}W^{q-1}(x)w(x)\,dx,
$$
that is, $\gamqv=\Gamma^q_1(W^{q-1}w)$.

\medskip
\noindent
(b) If $w(\infty)=0$, $1<q<\infty$ and $w_q$ is as in (\ref{ecuwq}), then
$\louw\subset\Gamma^q(w_q)\subset\gamuiw$. 

\medskip
\noindent
(c) If $w(\infty)>0$ then the embeddings $\louw\subset\gamqv\subset\gamuiw$ never hold.
\epro

\bdem
To prove (a) we use the  fact that  the embeddings $\louw\subset\gamqv\subset\gamuiw$ are equivalent to the
following condition (see
\cite{CPSS}):

\bequ\label{2inclu}
\bigg(\int_0^rw(x)\,dx\bigg)^q\approx\int_0^rv(x)\,dx+r^q\int_r^{\infty}v(x)\,{dx\over x^q},
\eequ
for every $r>0$. Hence,
$$
V(r)+r^q\int_r^{\infty}v(x)\,{dx\over x^q}\le C\int_0^rW^{q-1}(x)w(x)\,dx,
$$
which is equivalent to (see Theorem 2.2 in \cite{Ne})
$$
\int_0^{\infty}(f^{**}(x))^qv(x)\,dx\le C\int_0^{\infty}f^*(x)(f^{**}(x))^{q-1}W^{q-1}(x)w(x)\,dx.
$$
The other inequality follows similarly using now Theorem 3.2 in \cite{Ne}.

(b) follows immediately from Lemma~\ref{wq} and (\ref{2inclu}).

To finish, if we assume that the embeddings $\louw\subset\gamqv\subset\gamuiw$ hold, then by (\ref{2inclu}) and since
$w(\infty)=c>0$, we have that
$$
0<c^q\le C\liminf_{t\to\infty}{1\over t^q}\int_0^tv(x)\,dx,
$$
and hence $\int_t^{\infty}v(s)\,ds/s^q=\infty$, which is a contradiction.  $\qquad\qed$
\edem

We will also need to recall the following elementary observation.

\blem\label{meddec} 
If $\mu$ is a non-atomic positive measure in $[0,\infty)$ and $g$ is a positive decreasing function, then the mean value
function 
$$
{1\over \mu(0,t)}\int_0^tg(x)\,d\mu(x),
$$
is decreasing.
\elem

We are now ready to prove our main result.

\bteo\label{banachq1}
Let $1\le q<\infty$. The following are equivalent:

\medskip
\noindent
(a) $\gamquw$ is a Banach space.

\medskip
\noindent
(b)  $\gamquw\subset\gamqiw$.

\medskip
\noindent
(c) $w\in R_q$.
\eteo

\bdem
If $q=1$, this result is Theorem~\ref{normlouw}. Assume now that $1<q<\infty$. 
If $\Vert\cdot\Vert_{q,1}^*$ is a (rearrangement invariant) norm equivalent to $\Vert\cdot\Vert_{\gamquw}$, and since the
fundamental function $\Phi$ of $(\gamquw,\Vert\cdot\Vert_{q,1}^*)$ satisfies that
$$
\Phi(t)\approx W^{1/q}(t),
$$
then (see \cite{BS}):
$$
(\gamquw,\Vert\cdot\Vert_{\gamquw})=(\gamquw,\Vert\cdot\Vert_{q,1}^*)\subset \Gamma^{1,\infty}({\rm
d}\Phi)=\Gamma^{q,\infty}(w).
$$

Now if (b) holds, then by checking this embedding on characteristic functions we find that 
$$
\sup_{r<s}{W^{1/q}(s)\over s}\le C{W^{1/q}(r)\over r},
$$
and hence $w\in R_q$.  

To finish,  let us first observe that  $w\in R_q$ is equivalent to the condition $w/W^{1/q'}\in R_1$  and hence
$\Lambda^1(w/W^{1/q'})$ is Banach (Theorem~\ref{normlouw}). Thus, there exists a decreasing weight $\tilde w$ such that
$\Lambda^1(\tilde w)=\Lambda^1(w/W^{1/q'})$, and
$$
\int_0^r\tilde w(x)\,dx\approx \int_0^r{w(x)\over W^{1/q'}(x)}\,dx\approx\bigg(\int_0^rw(x)\,dx\bigg)^{1/q},
$$
and so,
$$
\int^r_0\widetilde W^{q-1}(x)\tilde w(x)\,dx\approx\bigg(\int_0^r\tilde w(x)\,dx\bigg)^q\approx \int_0^rw(x)\,dx.
$$
Thus, by Hardy's lemma (see \cite{BS}), 
$$
\int_0^{\infty}f^*(x)(f^{**}(x))^{q-1}\widetilde W^{q-1}(x)\tilde
w(x)\,dx\approx\int_0^{\infty}f^*(x)(f^{**}(x))^{q-1}w(x)\,dx,
$$
which proves that $\Gamma^q_1(w)=\Gamma^q_1(\widetilde W^{q-1}\tilde w)$. Therefore it suffices to study when
$\Gamma^q_1( W^{q-1} w)$ is Banach, assuming $w$ is decreasing.

If $w(\infty)=0$ then using Proposition~\ref{teo1} we find a weight $v$ such that $\Gamma^q_1( W^{q-1} w)=\gamqv$,
which is always a Banach space.

If  $w(\infty)=a>0$ and $w(0)=b<\infty$ then $w\approx 1$ and hence 
\begin{eqnarray*}
&&\bigg(\int_0^{\infty}f^*(x)(f^{**}(x))^{q-1}W^{q-1}(x)w(x)\,dx\bigg)^{1/q}\\
\approx&&\bigg(\int_0^{\infty}f^*(x)\bigg(\int_0^xf^*\bigg)^{q-1}\,dx\bigg)^{1/q}
\approx\int_0^{\infty}f^*(x)\,dx,
\end{eqnarray*}
and $\Gamma^q_1(W^{q-1} w)=L^1$.

If  $w(\infty)=a>0$ and $w(0)=\infty$, assuming without loss of generality that $a=1$, we define $\overline w=w-1$, and
$\overline W(t)=\int_0^t \overline w(x)\,dx$. Then 
\begin{eqnarray*}
&&\bigg(\int_0^{\infty}f^*(x)(f^{**}(x))^{q-1}W^{q-1}(x)w(x)\,dx\bigg)^{1/q}\\
\approx&&\bigg(\int_0^{\infty}f^*(x)(f^{**}(x))^{q-1}\Big[\overline W^{q-1}(x)\overline
w(x)+x^{q-1}\overline w(x)+\overline W^{q-1}(x)\Big]\,dx\bigg)^{1/q}\\
+&&
\bigg(\int_0^{\infty}f^*(x)(f^{**}(x))^{q-1}x^{q-1}\,dx\bigg)^{1/q},
\end{eqnarray*}
and hence,
$$
\Gamma^q_1( W^{q-1} w)=\Gamma^q_1(u)\cap L^1,
$$
where $u(x)=\overline W^{q-1}(x)\overline w(x)+x^{q-1}\overline w(x)+\overline W^{q-1}(x)$. Thus we only need to show
that $\Gamma^q_1(u)$ is Banach. We first of all observe that the function $u(t)/t^{q-1}$ is decreasing, and hence by
Lemma~\ref{meddec} we obtain that $t^{-q}\int_0^tu(x)\,dx$ is also a decreasing function. Thus, if we define
$V(r)=(\int_0^ru(x)\,dx)^{1/q}$, we have that $V$ is a positive increasing function such that $V(r)/r$  is decreasing, and
hence there exists a decreasing weight $\phi$ such that $V(r)\approx\Phi(r)=\int_0^r\phi(x)\,dx$.  Therefore, since
$$
\int_0^rV^{q-1}(x)v(x)\,dx\approx V^q(r)\approx\Phi^q(r)\approx\int_0^r\Phi^{q-1}(x)\phi(x)\,dx,
$$
we conclude that
$$
\Gamma^q_1(u)=\Gamma^q_1(V^{q-1}v)=\Gamma^q_1(\Phi^{q-1}\phi),
$$
and by the previous argument, if we show that $\phi(\infty)=0$ we are done. But, since $\phi$ is decreasing, then
\begin{eqnarray*}
\phi(\infty)&=&\lim_{t\to\infty}{\Phi(t)\over t}\le C\lim_{t\to\infty}{V(t)\over t}=
C\lim_{t\to\infty}\bigg({\int_0^tu(x)\,dx\over t^q}\bigg)^{1/q}\\
&\le&C\lim_{t\to\infty}\bigg({\int_0^t[\overline W^{q-1}(x)\overline w(x)+x^{q-1}\overline w(x)+\overline
W^{q-1}(x)]\,dx\over t^q}\bigg)^{1/q}\\
&\le& C\lim_{t\to\infty}\bigg({\overline W^{q}(t)\over t^q}+{\overline W(t)\over
t}+{\overline W^{q-1}(t)\over t^{q-1}}\bigg)^{1/q}=0,
\end{eqnarray*}
since $\overline w(\infty)=0$. $\qquad\qed$
\edem

\bobs\label{banachrp}
From the proof of the previous result, it is easy to show  that if $\gamqaw$ is a Banach space, then we always
have that $w\in R_q$.
\eobs

We consider next the normability characterization for the case $1<\alpha\le q$. We observe that now, we obtain the
same condition as for the $\Lambda^q(w)$ spaces.

\bteo\label{banachaq}
Let $1<\alpha\le q$. Then, $\gamqaw$ is a Banach space, if and only if, $w\in B_q$.
\eteo

\bdem
If $w\in B_q$ then by Theorem~\ref{normlop} we have that $\loqw=\gamqaw=\gamqw$. Conversely, if $\gamqaw$ is a
Banach space, then it is  a rearrangement invariant Banach function space, and hence (see \cite{BS})
$\gamqaw\subset\gamqiw$. In particular, if for  fixed $0<a<\infty$ and $s>1$, we choose the function
$$
f^*(x)=\cases{1 & if $0\le x\le a$,\cr$\smallbreak$
\displaystyle{a\over x} & if $a< x\le sa$,\cr
0 & if $sa<x$,}
$$
then this embedding gives
$$
f^{**}(x)W^{1/q}(x)\le C\bigg(\int_0^{\infty}(f^*(t))^{\alpha}(f^{**}(t))^{q-\alpha}w(t)\,dt\bigg)^{1/q},
$$
which for $x=as$ is equal to
\bequ\label{termiu}
\bigg({1+\log s\over s}\bigg)^qW(as)\le C\bigg(\int_0^aw(x)\,dx+\int_a^{as}\Big({a\over x}\Big)^{q}\Big(1+\log{x\over
a}\Big)^{q-\alpha}w(x)\,dx\bigg).
\eequ
On the other hand, if we choose the function $f^*=\chi_{[0,s]}$ we easily obtain that $w\in R_q$, and hence (Lemma~7.1 in
\cite{Ne}),
$$
\int_{a}^{as}\Big({a\over x}\Big)^{q}w(x)\,dx\le C(1+\log s)\int_0^aw(x)\,dx.
$$
If we combine this with (\ref{termiu}), we obtain,
$$
\bigg({1+\log s\over s}\bigg)^qW(as)\le C(1+\log s)^{q-\alpha+1}\int_0^aw(x)\,dx,
$$
i.e.,
$$
{W(as)\over (as)^q}\le C(1+\log s)^{1-\alpha}{W(a)\over a^q}.
$$
Finally, taking $s$ big enough so that $C(1+\log s)^{1-\alpha}<1$ (which we can do since $\alpha>1$), and using Lemma 6.3 in
\cite{Ne}, we obtain that $w\in B_q$.  $\qquad\qed$
\edem

We can also give a characterization for the case $q=1$ and $0<\alpha<1$.

\bteo\label{alphameu}
Let $w$ be a weight and $0<\alpha<1$. Then, $\gamuaw$ is a Banach space, if and only if, $w\in R_1$ and
$\gamuaw=\louw$.
\eteo

\bdem
If $\gamuaw$ is a Banach space, then by Remark~\ref{banachrp} we have that $w\in R_1$. Also, since the fundamental
function of $\gamuaw$ is like $W$, then $\louw\subset\gamuaw$, and the equality follows. The converse result is trivial by
Theorem~\ref{normlouw}. $\qquad\qed$
\edem

\bobs\label{runobancg}
Since $\Gamma^1_0(w)$ is always a Banach space, and $\Gamma^1_1(w)$ is Banach, if and only if $w\in R_1$, it would a
priori be natural to expect that the normability characterization for the space $\gamuaw$, $0<\alpha<1$, should be  weaker
than $R_1$, contrary to what we have just proved. In fact, it is easy to show that in general $w\in R_1$ does not imply
$\gamuaw=\louw$. For example, with $w=1$, and $f^*$ a bounded function, $f^*(t)=t^{-1}\log^{-1/\alpha}(t)$ for $t$ big
enough, one can prove that $f\in\louw\setminus\gamuaw$.
\eobs

\section{Equivalent norms and embeddings}

We consider in this section the existence of a norm in the space $\louw$, which depends on the maximal function $f^{**}$.
This question has a positive answer in all the other cases of normability: for the case  $\loqw$, $q>1$,  see
Theorem~\ref{normlop}, and for  $\loqiw$ see Theorem~\ref{normlopi}. Observe  that Theorem~\ref{normlouw} does
not give any information on whether there exists such an equivalent norm in $\louw$ (i.e., if it is a
Gamma space).  Since for $w$ decreasing, we have that
$\gamuw\subset\louw\subset\gamuiw$, and the endpoint spaces are normable in terms of the maximal function, we ask
ourselves when can we get equality on each case: 

\begin{enumerate}
\item[(i)] $\louw=\gamuw$, if and only if $w\in B_1$ (Theorem~\ref{bpam}).

\item[(ii)] $\louw=\gamuiw$,  if and only if, $w\in L^{\infty}$, and either $w(\infty)>0$ or, $w(\infty)=0$ and $w\in
L^1$ (see  \cite{CPSS}).
\end{enumerate}
It is now clear that there exist weights $w$ which do not satisfy either (i) or (ii). Thus, for the corresponding spaces
$\louw$, the problem of normability with respect to $f^{**}$ was an open question. The following result gives a positive
answer in all cases.

\bteo\label{gamuv}
Let $w$ be a decreasing weight. Then:
\medskip
 
\noindent
(i) There exists a weight $v$ such that $\louw=\gamuv$, if and only if $w(\infty)=0$.

\medskip
\noindent
(ii)  $\louw= L^1$,  if and only if $w(\infty)>0$ and $w\in L^{\infty}$.

\medskip
\noindent
(iii) There exists a weight $v$ such that $\louw=\gamuv\cap L^1$,  $L^1\nsubseteq \gamuv$, and $\gamuv\nsubseteq L^1$, 
if and only if
$w(\infty)>0$ and $w\notin L^{\infty}$.
\eteo

\bdem
Let us first prove (i).  $\louw=\gamuv$ is equivalent to the condition  (see \cite{CPSS}),

\bequ\label{eqwv}
{1\over r}\int_0^r w(s)\,ds\approx{1\over r}\int_0^r v(s)\,ds+\int^{\infty}_r {v(s)\over s}\,ds=S(S^*v)(r),
\eequ 
where $S^*$ is the adjoint operator of $S$:
$$
S^*f(r)=\int^{\infty}_r {f(s)\over s}\,ds.
$$
Suppose that $w(\infty)=\rho>0$.  Since $S(S^*v)$ is a decreasing function then,   using (\ref{eqwv}), we obtain the
existence of the limit 
$$
\lim_{r\to\infty}{1\over r}\int_0^rv(s)\,ds=\alpha\ge 0.
$$
Since $S(S^*v)=S^*(Sv)$, then we find that necessarily $\alpha=0$. But on the other hand,  
$$
\rho=\lim_{r\to\infty}{1\over r}\int_0^rw(s)\,ds\le C\lim_{r\to\infty}{1\over r}\int_0^rv(s)\,ds,
$$
and we get the contradiction that $\alpha\ge\rho/C>0$. Thus, $\rho=0$. 
\medskip

The converse is a direct consequence of (\ref{eqwv}) and Lemma~\ref{wq}.
\medskip

The proof of (ii) is trivial, since the condition says that $w\approx 1$. Let us finally consider (iii).  Observe that
$\louw\subset L^1$ is equivalent to $w(\infty)>0$, and hence, if we define $u(t)=w(t)-w(\infty)$, using (i) we find a weight
$v$ such that $\Lambda^1(u)=\gamuv$. Thus,  $\louw=\gamuv\cap L^1$. Also, if $L^1\subset\gamuv$,  then
$\louw= L^1$ which contradicts the fact that $w\notin L^{\infty}$.  Similarly, if  $\gamuv\subset L^1$ then (see
\cite{CPSS}) 
$$
t\le C\bigg(V(t)+t\int_t^{\infty}{v(s)\over s}\,ds\bigg),
$$
and hence, $v(\infty)>0$, which implies $\gamuv=\{0\}$. To finish, we only need to observe that if $\louw=\gamuv\cap L^1$
then, as we have mentioned before,  $w(\infty)>0$ and, if $w\in L^{\infty}$ then $w\approx 1$,  which says that
$\louw=L^1$, contradicting the hypothesis. $\qquad\qed$
\edem

\bobs
(i) Since $\Vert f\Vert_{L^1}=\sup_{t>0}(tf^{**}(t))$, it is now clear that Theorem~\ref{gamuv} shows that in $\louw$
there always exists an equivalent norm depending on the maximal function.

\medskip
(ii) For the weight $w=\chi_{(0,1)}$, one can show that the function obtained by the argument used in the previous proof
is $v(t)=t^{-2}(\log(4t)\chi_{(1/4,1/2)}(t)-\log t\chi_{(1/2,1)}(t))$. In fact, it is easy to prove that in order for
(\ref{eqwv}) to hold, it suffices that $v$ is a bounded function with compact support, which vanishes on a neighborhood of 0.
\eobs

We consider now other kind of embeddings for $\louw$.

\bpro\label{lugaq}
Let $w$ be a decreasing weight. Then the following are equivalent:

\medskip\noindent
(i) There exist $1<q<\infty$ and a weight $v$ such that $\louw=\gamqv.$

\medskip\noindent
(ii) For every $1<q<\infty$, $\louw=\Gamma^q_1(w_q)$, where $w_q$ is as in (\ref{ecuwq}).

\medskip
\noindent
(iii) $\louw=\gamuiw$, $w(\infty)=0$ and $w\in L^1$.
\epro

\bdem
It is easy to see that it suffices to show that (i) implies (iii). By Proposition~\ref{teo1} we have that
$\gamqv=\Gamma^q_1(W^{q-1}w)$, and hence
\begin{eqnarray*}
\int_0^{\infty}f^*(x)w(x)\,dx&\approx&\bigg(\int_0^{\infty}f^*(x)(f^{**}(x))^{q-1}W^{q-1}(x)w(x)\,dx\bigg)^{1/q}\\
&\le&\bigg(\bigg(\int_0^{\infty}f^*(x)w(x)\,dx\bigg)\sup_{t>0}\Big(f^{**}(t)W(t)\Big)^{q-1}\bigg)^{1/q},
\end{eqnarray*}
that is,
$$
\int_0^{\infty}f^*(x)w(x)\,dx\le \sup_{t>0}\Big(f^{**}(t)W(t)\Big),
$$
which shows that $\gamuiw\subset\louw$.  That $w(\infty)=0$, is also a consequence of Proposition~\ref{teo1}, and
the fact that $w\in L^1$ follows from \cite{CPSS}. $\qquad\qed$
\edem

Another kind of natural embeddings follow if we observe that on the one hand, if $w$ is decreasing,
$\gamuw\subset\louw\subset\gamuiw$, and on the other $\gamuw\subset\gamuqw\subset\gamuiw$, for $1<q<\infty$ (see
Definition~\ref{gamupw} (ii)). Thus, since $\gamuqw$ is a Banach space, we consider the question of whether
$\louw=\gamuqw$, for some $1<q<\infty$. The answer is rather surprising and says that this only happens in the trivial case
$w\in B_1$, i.e., when $\louw=\gamuw$.

\bteo\label{gamup}
Let $w$ be a decreasing weight, and $1<q<\infty$. Then, the following are equivalent:

\medskip\noindent
(i) $\louw\subset\gamuqw$.

\medskip\noindent
(ii) $S:\luwd\longrightarrow L^{1,q}(w)$.

\medskip\noindent
(iii) $w\in B_1$.
\eteo

\bdem
That (i) and (ii) are equivalent is easy. Also, if $w\in B_1$, then (see Theorem~\ref{bpam}) $S:\luwd\longrightarrow
L^{1}(w)\subset L^{1,q}(w)$, which proves (ii). Therefore it suffices to prove that (ii) implies (iii). A simple calculation shows
that if $w$ satisfies (ii), then $W^{q-1}w\in B_q$. Hence there exists $p<q$ such that $W^{q-1}w\in B_p$ (see \cite{AM}),
that is,
$$
\int_r^{\infty}{W^q(x)\over x^{p+1}}\,dx\le C{W^q(r)\over r^p},
$$
(this is the equivalent characterization of $B_p$ proved in \cite{So}).  Thus, for $s>1$,
$$
{W^q(sr)\over W^q(r)}\le Cs^p,
$$
and hence, if $s=2^k$, $k=1,2,\cdots,$
$$
\sup_{r>0}\bigg({W(2^kr)\over W(r)}\bigg)^{1/k}\le C^{1/k}2^{p/q},
$$
and using Theorems~4.2 and 4.5 in \cite{CGS} we prove that $w\in B_1$. $\qquad\qed$
\edem

A similar result is obtained if we study the equality between the spaces $\gamuqw$ and $\Gamma^q(w_q)$.

\bpro\label{gamqwgamqwq}
Let $w$ be a decreasing weight. Then, $w\in B_1$, if and only if, $\gamuqw=\Gamma^q(w_q)$, for every $1<q<\infty$.
\epro

\bdem
If $w\in B_1$ then $w(\infty)=0$, and by Proposition~\ref{teo1} we can find $w_q$ such that
$\louw\subset\Gamma^q(w_q)\subset\gamuiw$. Since $\louw\subset\gamuqw\subset\gamuiw$ also holds
(Theorem~\ref{gamup}), and
$\gamuqw=\Gamma^q(W^{q-1}w)$, then using again Proposition~\ref{teo1} we have that $\gamuqw=\Gamma^q(w_q)$. 
Conversely, 
\begin{eqnarray*}
&&\int_0^{\infty}(f^{**}(x))^qW^{q-1}(x)w(x)\,dx\approx\int_0^{\infty}f^*(x)(f^{**}(x))^{q-1}W^{q-1}(x)w(x)\,dx\\
&&\le\bigg(\int_0^{\infty}(f^{*}(x))^qW^{q-1}(x)w(x)\,dx\bigg)^{1/q}\bigg(\int_0^{\infty}(f^{**}(x))^{q}
W^{q-1}(x)w(x)\,dx\bigg)^{1/q'},
\end{eqnarray*}
and hence $\Lambda^q(W^{q-1}w)=\gamuqw$, which is a Banach space. Thus $W^{q-1}w\in B_q$, which is equivalent to $w\in
B_1$ (see, e.g., \cite{CM}). $\qquad\qed$
\edem

\bobs
If $w$ is a decreasing weight such that $w(\infty)=0$, then condition (\ref{ecuwq}) tells us that
$\gamuqw\subset\Gamma^q(w_q)$ (see \cite{CPSS}). Hence, as a consequence of Proposition~\ref{gamqwgamqwq}, we have
that under those hypotheses on the weight, $w\in R_1\setminus B_1$, if and only if,
$\gamuqw\varsubsetneq\Gamma^q(w_q)$.
\eobs

To finish, we give some examples which show how the space $\louw$ fits among  the range of Banach spaces
$\gamupw$.

\beje\label{ejemplo}
(i) Contrary to what happens in the case $\louw\subset\gamuqw$ (see Theorem~\ref{gamup}), the converse embedding can
hold true without having that necessarily  $\louw=\gamuiw$. For example, if $w(\infty)>0$ and $w(0)=\infty$, then  for
$1\le q<\infty$, 
$$
\{0\}=\gamuqw\varsubsetneq\louw\varsubsetneq\gamuiw.
$$

\medskip\noindent (ii) Another example can be obtained if we chose $w(t)=(1-\log t)\chi_{(0,1)}(t)$. In this case, if $1\le
q<\infty$, 
$$
\{0\}\neq\gamuqw\varsubsetneq\louw\varsubsetneq\gamuiw.
$$

\medskip\noindent (iii)  If $w(t)=\log t(\log t+2)\chi_{(0,e^{-2})}$, then, for $1<p<q<\infty$,
$$
\gamuw\varsubsetneq\gamupw\varsubsetneq\gamuqw\varsubsetneq\gamuiw,
$$
$\gamuw\varsubsetneq\louw\varsubsetneq\gamuiw$ and $\louw\nsubseteq\gamuqw$.

\medskip\noindent (iv) If $w(t)=t^{-\alpha}\chi_{(0,1)}(t)$, $-1<\alpha<0$, then, for $1<p<q<\infty$,
$$
\gamuw=\louw\varsubsetneq\gamupw\varsubsetneq\gamuqw\varsubsetneq\gamuiw.
$$
\eeje

\bigbreak

\end{document}